# The construction of all the projective spanning trees of a given digraph

Mikhail A. Antonets,  Grigoriy P. Kogan

In the present paper we consider the problem of constructing all the projective rooted spanning trees of a given graph. We propose an algorithm based on reducing this problem to the problem of constructing all the maximal independent sets of a certain derived graph.  We also offer an algorithm for step-by-step growing the spanning trees that sifts, at each step, any sub-trees non-extendable to projective spanning trees.  The proposed algorithms were applied for analyzing the possible variants of the syntactic subordination of nouns through prepositions  .

In the process of constructing the syntactic structure of a sentence according to Tesniere [1], the initial object is the digraph $G$ whose vertices are the word-forms of the sentence and whose arcs connect the pairs of word-forms that are able to have syntactic links not contradicting the grammar's rules. In such a process any projective spanning tree (see Definitions 1, 2) is a candidate to be selected as the syntactic structure of the sentence. However, even for not too long sentences the number of projective spanning trees is very large, while the time needed for finding them sharply grows with the sentence's length (see. [ 2-4]). Because of the latter fact, the construction of all the projective spanning trees of the digraph  $G$ might be actually only when the majority of the sentence's word-forms are already given the unique arcs entering them.

In the following Russian sentence, after determining all its syntactic links except those providing the connection of two nouns through a preposition, there exist several variants of choosing the subordination links through the four prepositions (i.e. going from a noun to a preposition or from a preposition to a noun) of this sentence, and only two of them seem reasonable:  «Ящик с опасным содержимым обнаружил участковый полиции во дворе жилого дома у грузчика из поселка Прудбой Калачевского района.» (However, due to some special rules of the English grammar, the above Russian sentence's translation into English doesn't suppose linguistic ambiguity).

Let us to start with some definitions and examples. Let $G$ be a graph with no loops whose vertices are indexed by the natural numbers from 1 to $n$, $A(G)$ be its adjacency matrix and $E$ be its edge set.

Let's correspond to every edge (or arc, if the graph is directed) $e$ of the graph $G$ a pair of natural numbers $\{\underline{m}, \overline{m}\}$, $\underline{m} < \overline{m}$, where $\underline{m}, \overline{m}$ are the indexes of the vertices connected by the edge. We'll call the numbers $\underline{m}, \overline{m}$ the beginning and the end of the edge $\{\underline{m}, \overline{m}\}$ correspondingly.

**Definition 1.** We'll say that the edges $\{\underline{m}, \overline{m}\}$ and $\{\underline{l}, \overline{l}\}$ are conflicted (or in conflict) if the intervals $[\underline{m}, \overline{m}]$, $[\underline{l}, \overline{l}]$ have a non-empty intersection different from any of both the intervals.

A geometrical illustration of the conflict of two edges is the fact of an edge intersection upon locating the vertices $\underline{m}, \overline{m}, \underline{l}, \overline{l}$ in the integer points (of a line) determined by their indexes and depicting the edges as strictly convex curves (or arcs) lying in the upper half-plane.

**Example.** Let's consider the digraph with vertices indexed by the natural numbers 1,…,7 that has the adjacency matrix

|   | 1 | 2 | 3 | 4 | 5 | 6 |
|---|---|---|---|---|---|---|
| 1 | 0 | 0 | 0 | 0 | 0 | 0 |
| 2 | 0 | 0 | 1 | 0 | 0 | 1 |
| 3 | 0 | 0 | 0 | 0 | 0 | 1 |
| 4 | 1 | 1 | 0 | 0 | 0 | 0 |
| 5 | 0 | 1 | 0 | 1 | 0 | 0 |
| 6 | 0 | 0 | 0 | 0 | 0 | 0 |

This graph can be depicted in the following way:

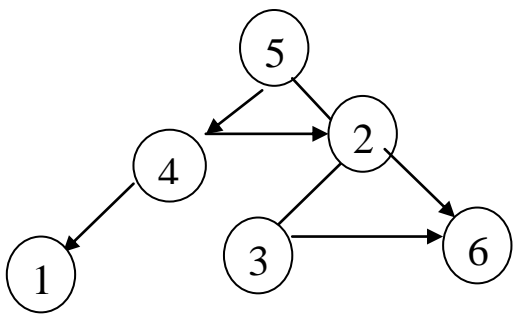

Fig. 1

The linear representation of the above graph has the form

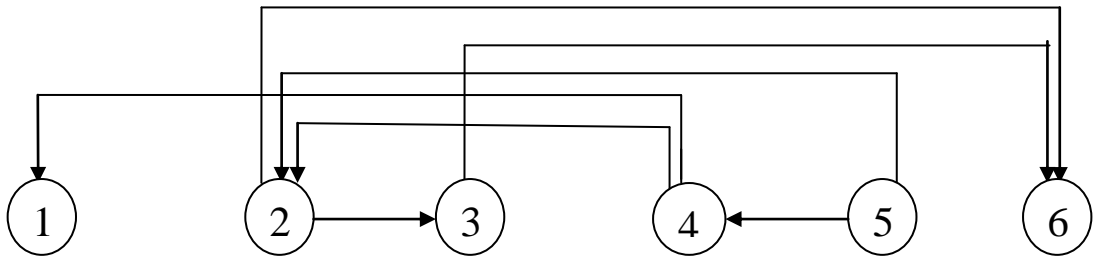

Fig. 2.

From the drawing (Fig. 2), one may conclude that there exist 5 points of edge intersection.

The set of conflicting arcs might be drawn as the following graph: to each conflicting arc, we correspond a vertex of the new graph indexed by the pair of the indexes of its beginning and end, and any two of such new vertices are connected by an edge if and only if the arcs corresponding to those new vertices are in conflict

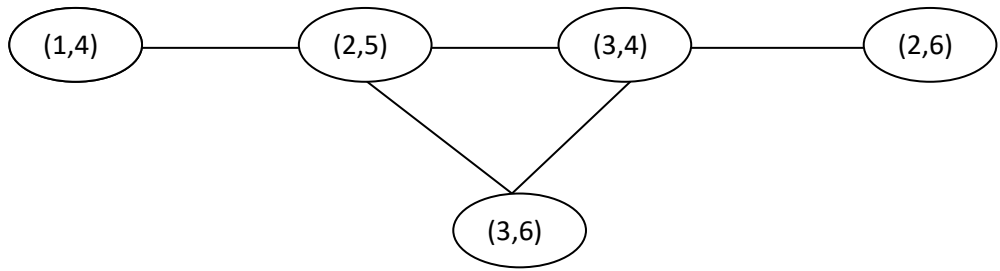

Fig.3.

**Definition 2.** A sub-graph $T$ of a graph $G$ will be called a maximal projective sub-graph of the graph $G$ if а) its set of vertices coincides with the vertex set of $G$;

б) any two edges of $T$ are not in conflict;

в) any edge of $G$ not belonging to $T$ is in conflict with at least one edge of $T$.

**Definition 3.** As the conflict graph $H$ of the graph $G$, we define the graph whose vertex set is $G$'s edge set $E$ an whose edge set is precisely all the pairs of conflicting edges of the graph $G$.

Let's denote the adjacency matrix of the above graph by $A(H)$.

It's easy to notice that a sub-graph $T$ is a maximal projective sub-graph in $G$ if and only if the edges of $T$ form a maximal independent sub-set of the vertex set of $H$.

**Problem. Find an efficient algorithm for the following problem:** given a digraph $\vec{G}$ with indexed vertices, construct all the directed projective spanning trees with a given root.

The following proposed algorithm supposes the construction, at the first stage, of all the maximal projective sub-graphs of the graph $\vec{G}$ and, afterwards, building, for each of those sub-graphs, all its directed spanning trees with the given root that, by the construction, should be projective.

**2.** In the further, we'll use the linear order on the vertex set of the conflict graph $H$ (i.e. on the edge set of the graph $G$) generated by the lexicographical order on this set:

$\{\underline{m}, \overline{m}\} > \{\underline{l}, \overline{l}\}$, if:

either $\overline{l} < \overline{m}$, or $\overline{l} = \overline{m}$ и $\underline{l} < \underline{m}$

This lexicographical order provides an opportunity to index the edges by the natural numbers k=1,..., $|E|$ via putting $k = K(\{\underline{m}, \overline{m}\})$ where the function $K$ on the edge set $E$ is given by the relationship

$$K(\{\underline{m}, \overline{m}\}) = \left|\{\{\underline{l}, \overline{l}\}: \{\underline{l}, \overline{l}\} \leq \{\underline{m}, \overline{m}\}\}\right| = \left|\{\{\underline{l}, \overline{l}\}: \overline{l} < \overline{m}\}\right| + \left|\{\{\underline{l}, \overline{m}\}: \underline{l} \leq \underline{m}\}\right|$$

The function $K$ strictly increases under the lexicographical order on the edge set $E$.

For a natural $k$, let's denote by $E^k$ the set of edges

$$E^k = \{\{\underline{m}, \overline{m}\}: \{\underline{m}, \overline{m}\} \in E, K(\underline{m}, \overline{m}) \leq k\} \qquad (1)$$

3. Let's notice that we can correspond to each sub-set of edges from $E$ a vector whose components are the indexes of edges in this set located in the order of increasing.

The following proposition provides a basis for describing the inductive process of constructing the set $L_{k+1}$ of all the maximal sub-sets of non-conflicting edges from the set $E^{k+1}$ by means of the set $L_k$ of all the maximal sub-sets of non-conflicting edges from the set $E^k$.

Let's denote by $B^{k+1}$ the set of edges from the set $E^k$ that are not conflicted with the edge $e_{k+1}$. Obviously, the edge $e_j$ belong to the set $B^{k+1}$ if and only if $A(H)_{r+1,j} = 0$.

For an arbitrary set $R$ of edges of the graph $G$, let's denote by $Vect(R)$ the binary vector from $\{0,1\}^{|E|}$ whose components are defined by the following relationship: $Vect(R)_j = 1$ if $e_j$ belongs to the set $R$ and $Vect(R)_j = 0$ otherwise.

**Proposition 1.** For an arbitrary maximal non-conflicted set $F$ of edges from the set $E^{k+1}$, there holds the following alternative:

1) $F = F'$, where $F'$ is a maximal non-conflicted set of edges from the set $E^k$ such that $F' \cap \overline{B^{k+1}} \neq \emptyset$;
2) there exists at least one maximal non-conflicted set $F'$ of edges from the set $E^k$ such that for any $l$, such that $e_l \in B^{k+1} \setminus F'$, there holds the inequality
$$\sum_{e_j \in F' \cap B^{k+1}} A(H)_{l,j} \geq 1 \qquad (2)$$

**Proof.** For the structure of an arbitrary maximal non-conflicted set $F$ of edges from $L_{k+1}$, only the following two options are possible:

1) either $F = F'$, where $F'$ belongs to $L_k$ and contains at least one edge conflicted with the edge $e_{k+1}$, while it's possible if and only if the condition (1) of the considered proposition is true;
2) or $F = (\overline{F} \cap B^{k+1}) \cup \{e_{k+1}\}$, where $\overline{F} \in L_k$ is such that the intersection $\overline{F} \cap B^{k+1}$ is one of the maximal-by-inclusion subsets of the set $\{F' \cap B^{k+1}: F' \in L_k\}$.

But a sub-set of the form $F' \cap B^{k+1}$ is maximal-by-inclusion in the set $\{F' \cap B^{k+1}: F' \in L_k\}$ if and only if for any edge $e_j$ from the set $B^{k+1} \setminus F'$ there exists at least one conflicted (with that edge) edge $e_l$ from the set $F' \cap B^{k+1}$, i.e. the entry $A(H)_{lj}$ of the adjacency matrix of the graph $H$ equals 1, and, therefore, the condition (2) holds. The proposition is proved.

## An algorithm for constructing all the maximal sub-sets of non-conflicting edges of the graph $G$.

Let $e_k$ be the edge with index $k$ and $L_k = \{F^k_1, \ldots, F^k_{l_k}\}$ be the set of all the maximal sub-sets of non-conflicting edges from the set $E^k$.

Let's construct the set of all the maximal sets of non-conflicting edges from the set $E^{k+1}$.

For $k = 1$ we get the equality $F^1_1 = E^1$ and the set of sub-sets $L_k$ has the form

$$L_1 = \{F^1_1\}$$

where $F^1_1 = \{e_1\}$, $e_1$ is the youngest edge of the graph.

Now let's assume that the set $L_k = \{F^k_1, \ldots, F^k_{l_k}\}$ of all the maximal sub-sets of non-conflicting edges from the $E^k$ was already built.

We'll construct the set $L_{k+1} = \{F^{k+1}_1, \ldots, F^{k+1}_{l_{k+1}}\}$ via traversing the sets from $L_k$.

Let $F^k_j$ be a current set from $L_k$.

If the inclusion $F^k_j \subseteq B^{k+1}$ isn't true then, due to statement (1) of the Proposition 1, the set $F^k_j$ is an element of the set $L_{k+1}$;

and, if the condition (2) holds for $F' = F^k_j$, to the set $L_{k+1}$ we add the set of edges $(F^k_j \cap B^{k+1}) \cup \{e_{k+1}\}$, if it wasn't yet added to this set earlier. We repeat the described procedure for each set $F^k_j$ from the set $L_k$. The end of the algorithm.

**Comment.** It seems pretty reasonable to number, at the end of each stage of the above algorithm, all the constructed sets of edges in the accordance with the lexicographical order.

Because each projective spanning tree of the graph $G$ is a spanning tree of some maximal projective sub-graph of the graph $G$, while finding all the projective spanning trees of all the constructed maximal projective sub-graphs of the graph $G$, we'll find all its projective spanning trees both in the directed and undirected cases.

## The algorithm of construction of all the projective spanning trees

Let's denote by $A^{in}$ the in-degree matrix of graph G defined by the equation

$$A^{in} = -A(G) + \Lambda^{in}$$

where $A(G)$ is G's adjacency matrix, the matrix $\Lambda^{in}$ is defined as

$$\Lambda_{ij}^{in} = \frac{1}{2} deg_i^{in} \delta_{ij},$$

and $deg_i^{in}$ is the number of arcs entering the vertex $i$. This number equals the number of $A$'s entries in its $i$-th row that are equal to 1.

For the number $NST(r)$ of G's spanning trees having the vertex $r$ as their root, we have the following formula:

$$NST(r) = a_r^{in}$$

where $a_r^{in}$ is $A^{in}$'s algebraic complement of its entry $A_{rr}^{in}$.

Finding all the projective spanning trees of a digraph with a given adjacency matrix and a given root.

Let us given a weighted digraph $G = (V, E)$ with n vertices indexed by integer numbers where a vertex $r$ is selected to be called the root.

**The task:** to find an algorithm for constructing all the marked projective spanning trees of G.

**The solution:** by induction, let's define the sets $T(k, r)$ of G's projective pre-spanning trees of generation k. By the definition, the set $T(0, r)$ of $G$'s projective pre-spanning trees consists of the root $r$ if the digraph $G$ has at least one $r$–rooted spanning tree.

The program for verifying the latter condition is Solver G.1.

Let $T(k,r)$ be the set $r$-rooted projective pre-spanning trees of generation $k$ of the graph $G$. Then, by the definition, the set $T(k+1,r)$ of $G$'s $r$-rooted projective pre-spanning trees of generation k+1 consists of all the $r$-rooted projective spanning trees of G that satisfy the following conditions:

1) the lengths of those trees are $k+1$;

2) any tree $t'$ from the set $T(k+1,r)$ can be received from some projective pre-spanning trees $t$ from the set $T(k,r)$ via adding new arcs going from the latter trees' leafs whose distance to the root $r$ is $k$ and having no conflicts either between each other or with the arcs of the pre-spanning tree $t$;

3) each tree $t'$ from the set $T(k+1,r)$ can be extended to a spanning tree of the graph $G$ via connecting its leafs whose distance to the root $r$ is $k+1$ with the roots of a forest of the sub-graph $G \setminus t'$.

## The algorithm for verifying the fulfillment of the condition 3) for a tree $t'$ of length $k+1$:

let's construct a new graph $G(t')$ received from the graph $G \setminus t'$ (i.e. the sub-graph of $G$ induced by the set of vertices $V(G) \setminus V(t')$) via adding a new vertex $v^*$ and arcs going from $v^*$ and entering each vertex of $G \setminus t'$ entered by at least one arc going from a leaf of the tree $t'$ whose distance from the root $r$ is equal to $k+1$.

Hence the graph $G(t')$ is defined by its set of vertices

$$V(G(t')) = (V(G) \setminus V(t')) \cup \{v^*\}$$

and its set of arcs

$$E(G(t')) = E(G \setminus t') \cup E_k,$$

where

$$E_k = \{(v^*, z), z \in Z_k\}, Z_k = \{z: \exists d \in V(t'), dist(r,d) = k+1, (d,z) \in E(G)\}$$

Let $A(t')$ be the adjacency matrix of the graph $G(t')$ and the matrix $\Lambda(t')^{in}$ be defined as

$$\Lambda(t')_{ij}^{in} = \frac{1}{2} \deg(t')_i^{in} \delta_{ij},$$

where $\deg(t')_i^{in}$ is the number of arcs entering the vertex $i$ (this number equals the number of entries of the matrix $A(t')$ in its $i$-th row that are equal to 1).

First we should verify the existence of a directed $v^*$-rooted spanning tree in $G(t')$ via computing the algebraic complement of the $v^*, v^*$-entry of its in-degree matrix $A(t')^{in}$:

$$A(t')^{in} = -A(t') + \Lambda(t')^{in}$$

Hence, for the number $NST(v^*, t')$ of spanning trees of the graph $G(t')$ having the vertex $v^*$ as their root, we have the following formula:

$$NST(v^*) = a(t')_{v^*}^{in}$$

where $a(t')_{v^*}^{in}$ is $A(t')^{in}$'s algebraic complement of its entry $A(t')_{v^*v^*}^{in}$.

This value is equal to the number of spanning trees of the graph G, containing the tree $t'$. The proof of this statement is based on the well-known analogue of Kirchhoff's theorem for directed graphs that states that the number of spanning trees of a given digraph rooted in the vertex i equals the determinant of its in-degree matrix with the i-th row and column removed.

(see http://stu.alnam.ru/book_grnet-53)

The tree $t'$ belongs to $T(k + 1, r)$ (i.e. is a part of some spanning tree) if and only if the above-mentioned determinant isn't zero.

References.